\documentclass{amsart}
\usepackage{amssymb, latexsym, amsmath}

\usepackage[all]{xy}

\usepackage[usenames, dvipsnames]{color}

\usepackage{mathrsfs} 

\usepackage[OT2,T1]{fontenc}
\DeclareSymbolFont{cyrletters}{OT2}{wncyr}{m}{n}
\DeclareMathSymbol{\Sha}{\mathalpha}{cyrletters}{"58}

\theoremstyle{plain}
\newtheorem{Theorem}[equation]{Theorem}
\newtheorem{Lemma}[equation]{Lemma}
\newtheorem{Corollary}[equation]{Corollary}

\newtheorem*{Conjecture*}{Conjecture}

\theoremstyle{definition}

\newtheorem{Paragraph}[equation]{}

\theoremstyle{remark}

\numberwithin{equation}{section}

\renewcommand{\k}{{\kappa}}

\newcommand{\D}{{\mathscr D}}

\newcommand{\F}{{\mathbb F}}

\renewcommand{\H}{{\text{\rm H}}}

\newcommand{\lL}{{{}_\ell\text{\rm L}}}
\newcommand{\nL}{{{}_n\text{\rm L}}}

\renewcommand{\O}{{\text{\rm O}}}

\newcommand{\Q}{{\mathbb Q}}

\renewcommand{\S}{{\mathcal S}}

\newcommand{\Z}{{\mathbb Z}}

\newcommand{\br}[1]{{\left<{#1}\right>}}

\newcommand{\car}{{\text{\rm char}}}

\newcommand{\df}{{\,\overset{\text{\rm df}}{=}\,}}

\renewcommand{\div}{{\text{\rm div}}}

\newcommand{\ind}{{\text{\rm ind}}}
\renewcommand{\inf}{{\text{\rm inf}\,}}

\newcommand{\isom}{{\;\simeq\;}}

\newcommand{\per}{{\text{\rm per}}}

\newcommand{\red}{{\text{\rm red}}}

\newcommand{\res}{{\text{\rm res}}}

\newcommand{\Br}{{\text{\rm Br}}}

\newcommand{\nBrdim}{{{}_n\text{\rm Br.dim}}}
\newcommand{\lBrdim}{{{}_\ell\text{\rm Br.dim}}}

\newcommand{\Spec}{{\text{\rm Spec}\,}}

\begin{document}

\title[Tame division algebras over
function fields of $p$-adic curves]
{Tame division algebras of prime period over
function fields of $p$-adic curves}


\author{Eric Brussel and Eduardo Tengan}
\address 
{Department of Mathematics \& Computer Science\\
Emory University\\
Atlanta, GA 30322\\ USA}
\email{brussel@mathcs.emory.edu}
\address
{Instituto de Ci\^encias Matem\'aticas e de Computa\c c\~ao\\
Universidade de S\~ao Paulo\\
S\~ao Carlos, S\~ao Paulo\\ Brazil}
\email
{etengan@icmc.usp.br}





\begin{abstract}
Let F be a field finitely generated and of transcendence degree one over a $p$-adic field, and 
let $\ell\neq p$ be a prime.  Results of Merkurjev and Saltman
show that $\H^2(F,\mu_\ell)$ is generated by
$\Z/\ell$-cyclic classes.  We prove the ``$\Z/\ell$-length'' in $\H^2(F,\mu_\ell)$
equals the $\ell$-Brauer dimension, which Saltman showed to be two.
It follows that all $F$-division algebras of period $\ell$ are crossed products, either
cyclic (by Saltman's cyclicity result) or tensor products of two cyclic division algebras.
Our result was originally proved by Suresh when $F$ contains the $\ell$-th roots of unity $\mu_\ell$.
\end{abstract}

\maketitle

\section{Introduction}

Let $F$ be a field and $n$ a number prime-to-$\car(F)$.
Suppose the cup product map
$\H^1(F,\mu_n)\otimes_\Z\H^1(F,\Z/n)\to\H^2(F,\mu_n)={}_n\Br(F)$ is surjective.
This is the case when $F$ contains the $n$-th roots of unity or $n=3$, by Merkurjev-Suslin's
theorem \cite[Theorem 16.1, Corollary 16.4]{MS83}, but unknown in general.
Let the {\it $n$-Brauer dimension} $\nBrdim(F)$ denote the smallest number $d$
such that every class in $\H^2(F,\mu_n)$ has index dividing $n^d$, 
and let the {\it $\Z/n$-length} $\nL(F)$ denote the smallest
number of $\Z/n$-cyclic classes needed to write any class in $\H^2(F,\mu_n)$.
It is easy to see that $\nBrdim(F)\leq\nL(F)$, but even when the cup product
map is surjective it is not known
in general whether finite $n$-Brauer dimension implies finite $\Z/n$-length.

Suppose $F$ is finitely generated of transcendence degree one over the $p$-adic field $\Q_p$.
Saltman showed $\nBrdim(F)=2$ in \cite[Theorem 3.4]{Sa97}.
Recently Suresh showed that $\nL(F)=2$ when $n=\ell$ is prime and $F$ contains $\mu_\ell$
(\cite[Theorem 2.4]{Sur10}).  
The assumption on roots of unity excludes important cases such
as the rational function field $F=\Q_p(T)$ (if $\ell\neq 2,3$ and $p\neq 1\pmod\ell$).
But in this case the cup product map is surjective by Merkurjev's theorem 
\cite[Theorem 2]{Mer83} and Saltman's cyclicity result for classes of prime index \cite[Theorem 5.1]{Sa07},
so $\H^2(F,\mu_\ell)$ is generated by $\Z/\ell$-cyclic classes.
We show $\lL(F)=2$, 
hence all $F$-division algebras of period $\ell$ and index $\ell^2$
decompose into two cyclic $F$-division algebras of index $\ell$.
It follows immediately that all $F$-division algebras of period $\ell$ are (abelian) crossed products.
Noncrossed products of larger period exist by \cite{BMT} and \cite{BT11}.

Our results rely on Saltman's degree-$\ell$ cyclicity result  
and hot point criterion \cite[Corollary 5.2]{Sa07}, and also on our lifting results
from \cite{BT11}, which use the machinery of Grothendieck's existence theorem.
We show an $F$-division algebra $\Delta$ of period $\ell$ and index $\ell^2$ is decomposable
by constructing a tensor factor of degree $\ell$, lifting a class constructed over
the generic points of the closed fiber of a 2-dimensional model $X/\Z_p$,
as in \cite{BT11}.  
We use Grunwald-Wang's theorem
to construct the lifted class so that it ``cancels'' the hot points of $\Delta$, which implies it is
part of a decomposition of $\Delta$ by the hot point criterion.
The cyclicity result then shows the remaining factor is cyclic.
Suresh's approach in \cite{Sur10} similarly cancels $\Delta$'s hot points using a tensor factor,
but his tensor factor is constructed
as a symbol algebra, which requires $\mu_\ell\subset F$.
We know of no way to get to the general case from that construction.

\section{Background and Conventions}

\Paragraph{Brauer Group Conventions.}\label{brauer}
In this paper an {\it $F$-division algebra} is a division ring that is central and finite-dimensional
over $F$.  If $D$ is an $F$-division algebra we write $[D]$ for the class of $D$ in 
the Brauer group $\Br(F)$,
$\ind(D)$ for the {\it index} or {\it degree} of $D$, and $\per(D)$ for the {\it period} of $D$.
We say $D$ is a {\it crossed product} if it contain a maximal subfield that is Galois over $F$.
See \cite{ABGV} for a discussion of crossed product and noncrossed product division algebras.

We write $\H^2(F,\mu_n)={}_n\Br(F)$ for the $n$-torsion
subgroup, where $n$ is prime-to-$\car(F)$ and $\mu_n$ is the group of $n$-th roots of unity.
In the terminology of \cite[Section 4]{ABGV}, the {\it $n$-Brauer dimension} $\nBrdim(F)$ 
is the smallest number $d$ such that every class in $\H^2(F,\mu_n)$ has index dividing $n^d$.  
We say $\H^2(F,\mu_n)$ is {\it generated by $\Z/n$-cyclic classes} if the cup product
map $\H^1(F,\mu_n)\otimes_\Z\H^1(F,\Z/n)\to\H^2(F,\mu_n)$ is surjective, and a class
is {\it $\Z/n$-cyclic} if it has the form $(f)\cdot\theta$ for some $(f)\in\H^1(F,\mu_n)$
and $\theta\in\H^1(F,\Z/n)$.
When $\H^2(F,\mu_n)$ is generated by $\Z/n$-cyclic classes,
we define the {\it $\Z/n$-length} $\nL(F)$ to be the smallest number $d$ such that
every element of $\H^2(F,\mu_n)$ can be written as a sum of $d$ $\Z/n$-cyclic classes.
See \cite[Section 3]{ABGV} for a discussion of known results regarding $\Z/n$-length,
usually called ``symbol length'' when $F$ contains an $n$-th root of unity.

If $\nL(F)=d$
then it is clear that $\nBrdim(F)\leq d$.  Conversely, it is not known
whether a finite Brauer dimension implies a finite $\Z/n$-length, or even that 
$\H^2(F,\mu_n)$ is generated by $\Z/n$-cyclic classes.
However, when $n=\ell$ Merkurjev proved that $\H^2(F,\mu_\ell)$ is generated by classes of
index $\ell$ (\cite[Theorem 2]{Mer83}), 
hence for the fields considered in this paper $\H^2(F,\mu_\ell)$ is
generated by $\Z/\ell$-cyclic classes by Saltman's cyclicity result.

\Paragraph{General Conventions.}
Let $S$ be an excellent scheme and suppose $n$ is invertible on $S$.
We write $\Z/n(r)$ for the \'etale sheaf $\Z/n$ twisted by an integer $r$,
and $\H^q(S,r)=\H^q(S,\Z/n(r))$ for the \'etale cohomology group.  
If $S=\Spec A$ for a ring $A$, we write $\H^q(A,r)$.
If $T$ is a subscheme of $S$ we write $\kappa(T)$ for its ring of meromorphic functions, which
is the localization of $\O_T$ at its associated points.
If $T\to S$ is a morphism of schemes then the restriction $\res_{S|T}:\H^q(S,r)\to\H^q(T,r)$
is defined, and we write $\beta_T=\res_{S|T}(\beta)$, and if $T=\Spec B$ we write $\beta_B$.
If $Z\subset S$ is another morphism we write $Z_T$ for the fiber product $Z\times_S T$.

If $v$ is a valuation on a field $F$, we write $\kappa(v)$ for the residue field of the
valuation ring $\O_v$, and $F_v$ for the completion of $F$ with respect to $v$.
If $v$ arises from a prime divisor $D$ on $S$, we write $v=v_D$, $\kappa(D)$, and $F_D$.
If a set $\{v_i\}$ arises from a divisor $D=\sum_i D_i$, we write $F_D=\prod_i F_{D_i}$.
Recall that if $F=(F,v)$ is a discretely valued field and $\alpha\in\H^q(F,r)$, 
then $\alpha$ has a {\it residue} $\partial_v(\alpha)$ in $\H^{q-1}(\kappa(v),r-1)$.
More generally if $\xi$ is a generic point of a scheme $S$, $F=\kappa(\xi)$, and
$\alpha\in\H^q(S,r)$, then for each discrete valuation $v$ on $F$ we define
\[\partial_v(\alpha)\df\partial_v(\alpha_F)\in\H^{q-1}(\kappa(v),r-1)\] 
We say $\alpha$ is {\it unramified} with respect to $v$ if $\partial_v(\alpha)=0$,
and in that case the {\it value} of $\alpha$ at $v$ is
the element $\alpha(v)=\res_{F|F_{v}}(\alpha)\in\H^q(\kappa(v),r)\leq\H^q(F_v,r)$
(\cite[7.13, p.19]{GMS}).
If $v$ arises from a prime divisor $D$ on a scheme, we will substitute the notations
$\partial_D$ and $\alpha(D)$.
If $S$ is noetherian
we write $D_\alpha$ for the {\it ramification divisor of $\alpha$ on $S$}, which is the sum of
(finitely many) prime divisors on $S$ at which $\alpha$ ramifies.

\Paragraph{Setup.}\label{setup}
In the following, 
$F$ will always be a finitely generated field extension of $\Q_p$ of transcendence degree one,
$n$ will be a prime-to-$p$ number,
and $X/\Z_p$ will be a connected regular (projective, flat) 
relative curve over $\Z_p$ with function field $F=K(X)$.
Such a surface exists for any $F$ by a theorem of Lipman (see \cite[Theorem 8.3.44]{Liu}).
We write $X_0=X\otimes_{\Z_p}\F_p$ for the closed fiber,
$C=X_{0,\red}$ for the reduced scheme underlying the closed fiber,
$C_1,\dots,C_m$ for the irreducible components of $C$, and $\S$ for the set
of singular points of $C$.
We assume that $X_0$ has normal crossings, hence that each $C_i$ is regular, 
and at most two of them meet (transversally) at each singular point of $C$.
This is situation is always achievable, by embedded resolution of curves in surfaces 
(see \cite[Theorem 9.2.26]{Liu}).

We say an effective divisor $D$ on $X$ is {\it horizontal} 
if each of its irreducible components maps surjectively to $\Spec\Z_p$.
By \cite[Proposition 2.6]{BT11} there exists for each closed point $z\in X\backslash\S$
a regular irreducible horizontal divisor $D\subset X$ that intersects $C$ transversally.
Let $\D_\S$ denote the support of these lifts.
We say a divisor $D$ is {\it distinguished} and write $D\in\D_\S$ if it is reduced
and supported in $\D_\S$.
Each $D\in\D_\S$ is a disjoint union of its irreducible components, each of which
has a single closed point and meets $C$ transversally.

By weak approximation (\cite[Lemma]{Sa98})
there exists an element $\pi\in F$ such that $\div(\pi)=C+E\subset X$,
where $E$ is horizontal and avoids all closed points of any finite set.
We fix such a $\pi$ such that $\div(\pi)$ avoids $\S$.

Suppose $\ell\neq p$ is prime, $\alpha\in\H^2(F,\mu_\ell)$, and $D_\alpha\subset X$ has normal crossings.
Following Saltman's terminology in \cite{Sa07} we say $\alpha$ has a {\it hot point} $z$ on $X$
if (and only if) $z$ is a nodal point of $D_\alpha$, and if $D,D'\subset D_\alpha$
are the two irreducible components meeting transversally at $z$,
then $\partial_D(\alpha)$ and $\partial_{D'}(\alpha)$ are unramified at $z$, and
$\br{\partial_D(\alpha)(z)}\neq\br{\partial_{D'}(\alpha)(z)}$.
By \cite[Corollary 5.2]{Sa07}, $\alpha$ has index $\ell$ if and only if $D_\alpha$ 
has no hot points (hot point criterion), 
and by \cite[Theorem 5.1]{Sa07}, if $\alpha$ has index $\ell$ then it is cyclic.

\Theorem[{\cite[Lemma 4.6, Theorem 4.9]{BT11}}]\label{lambda}
Assume the setup of \eqref{setup}.

{\rm a)} 
The image of the natural map $\H^1(\O_{C,\S},\Z/n)\to\H^1(\k(C),\Z/n)$ consists of  
the set of tuples $(\theta_1,\dots,\theta_m)\in\H^1(\kappa(C),\Z/n)$
such that each $\theta_i$ is unramified at each
$z\in\S\cap C_i$, and $\theta_i(z)=\theta_j(z)\in\H^1(\kappa(z),\Z/n)$ whenever $z\in C_i\cap C_j$.

{\rm b)}  For $q\geq 0$ and any integer $r$
there is a map $$\lambda:\H^q(\O_{C,\S},r)\to\H^q(F,r)$$ 
and a commutative diagram
\[
\xymatrix{
\H^q(\O_{C,\S},r)\ar[r]^-\lambda\ar[d]_\res&\H^q(F,r)\ar[d]^\res\\
\H^q(\kappa(C),r)\ar[r]^-\inf&\H^q(F_C,r)
}\]
such that if $\alpha_C\in\H^q(\O_{C,\S},r)$ and $\alpha=\lambda(\alpha_C)$
then:
\begin{enumerate}
\item[i)]
$\alpha$ is defined at the generic points of $C_i$, and $\alpha(C_i)=\res_{\O_{C,\S}|\k(C_i)}(\alpha_C)$.
\item[ii)]
The ramification locus of $\alpha$ (on $X$) is contained in $\D_\S$.
\item[iii)]
If $D\in\D_\S$ is prime and $z=D\cap C$,
then $\partial_D\cdot\lambda=\inf_{\kappa(z)|\kappa(D)}\cdot\partial_z$.
\item[iv)]
If $\alpha_C$ is unramified at a closed point $z$,
and $D$ is any (horizontal) prime lying over $z$,
then $\alpha$ is unramified at $D$, and has value
$\alpha(D)=\inf_{\kappa(z)|\kappa(D)}(\alpha_C(z))$.
\end{enumerate}
\rm

\section{Computations}

We first construct the cyclic class $\gamma\in\H^2(F,\mu_n)$ using a lift from $\H^1(\kappa(C),\Z/n)$.

\Lemma\label{lemma}
Assume the setup of \eqref{setup}.
Suppose $\theta_C\in\H^1(\O_{C,\S},\Z/n)$ maps to $(\theta_1,\dots,\theta_m)\in\H^1(\kappa(C),\Z/n)$
as in Theorem~\ref{lambda}(a), 
such that $\theta_C$ is unramified at all $z\in E\cap C$, with value $\theta_C(z)=0$
(in addition to being unramified at $\S$).
Let $\gamma=(\pi)\cdot\lambda(\theta_C)\in\H^2(F,\mu_n)$.  
Then for a prime divisor $D$ on $X$
$$
\partial_D(\gamma)=\begin{cases}
\theta_i &\text{ if $D=C_i$}\\
-\inf_{\kappa(z)|\kappa(D)}(\partial_{z}(\theta_C))\cdot(\pi)
&\text{ if $D\in\D_\S$, $z=D\cap C$, and $v_D(\pi)=0$.}\\
0&\text{ otherwise}\\
\end{cases}
$$
The ramification divisor $D_\gamma$ has normal crossings,
and consists of each $C_i$ at which $\theta_i$ is nonzero,
together with all $D\in\D_\S$ lifting $z:\partial_z(\theta_C)\neq 0$.
\rm

\begin{proof}
Set $\theta=\lambda(\theta_C)$.
Let $D\subset X$ be a prime divisor and let $z=D\cap C$.
Then
$$
\partial_D(\gamma)=\big[v_D(\pi)\theta-\partial_D(\theta)\cdot(\pi)
+v_D(\pi)\partial_D(\theta)\cdot(-1)\big]_{F_D}
$$
This element is in the subgroup $\H^1(\kappa(D),\Z/n)\leq\H^1(F_D,\Z/n)$.  

If $D=C_i$ is an irreducible component of $C$ then since $E$ contains no components of $C$
(by \eqref{setup}) we
have $v_D(\pi)=1$, and since $\theta_C\in\H^1(\O_{C,\S},\Z/n)$ we have
$\partial_D(\theta)=0$ by Theorem~\ref{lambda}(b)(ii).
Therefore $\partial_D(\gamma)=\res_{F|F_{C_i}}(\theta)=\theta_i\in\H^1(\kappa(C_i),\Z/n)$.
If $D\in\D_\S$ and $v_D(\pi)=0$
then $\partial_D(\gamma)=-\partial_D(\theta)\cdot(\pi)=-\inf_{\kappa(z)|\kappa(D)}(\partial_{z}(\theta_C))\cdot(\pi)$
by Theorem~\ref{lambda}(b)(iii).

It remains to show $\partial_D(\gamma)=0$ if $D$ is horizontal and runs through a point of $\S$,
if $D$ is horizontal and $v_D(\pi)\neq 0$, or if $D\not\in\D_\S$.
If $D$ is horizontal and runs through a point of $\S$,
then $v_D(\pi)=0$ since $E$ avoids $\S$, and $\partial_D(\theta)=0$ by Theorem~\ref{lambda}(b)(ii),
hence $\partial_D(\gamma)=0$.
If $D$ is horizontal, avoids $\S$, 
and $v_D(\pi)\neq 0$, then $D$ is a component of $E$,
so by assumption, $\partial_z(\theta_C)=0$ for $z\in E\cap C$ and $\theta_C(z)=0$.
Thus $\partial_D(\theta)=0$ and
$\theta(D)=0$ by Theorem~\ref{lambda}(b)(iv).
Therefore $\partial_D(\gamma)=0$.
Finally, if $D$ is horizontal, avoids $\S$,
$v_D(\pi)=0$, and $D\not\in\D_\S$ then $\partial_D(\gamma)=-\partial_D(\theta)\cdot(\pi)=0$
since $\partial_D(\theta)=0$ by Theorem~\ref{lambda}(b)(ii).

Suppose $D\in\D_\S$ and $v_D(\pi)=0$,
so $\partial_D(\gamma)=-\inf_{\kappa(z)|\kappa(D)}(\partial_{z}(\theta_C))\cdot(\pi)$, where $z=D\cap C$.
Since $v_D(\pi)=0$, $\pi$ is a local equation for $C$ at $z$
and $D$ intersects $C$ transversally at $z$, the image of $\pi$ in the local field $\kappa(D)$ is a uniformizer, 
hence $(\pi)$ has order $n$ in $\H^1(\kappa(D),\mu_n)$.
Thus $\partial_D(\gamma)$ is nonzero in this case if and only if $\partial_z(\theta_C)\neq 0$.

We conclude $D_\gamma$ consists of the components 
$C_i$ of $C$ for which $\theta_i$ is nonzero, together with the distinguished 
prime divisors $D\in\D_\S$ lying over points $z$ at which $\theta_C$ is ramified.
Since all such $D$ are regular and intersect $C$ transversally, $D_\gamma$ has normal crossings.
\end{proof}

Next we show an $F$-division algebra of prime period $\ell\neq p$ and index $\ell^2$
is decomposable
by constructing a cyclic factor using Lemma~\ref{lemma}, designed to cancel the division algebra's 
hot points.

\Theorem\label{theorem}
Let $F$ be a field finitely generated of transcendence degree one over $\Q_p$, and
suppose $\Delta$ is an $F$-division algebra of prime period $\ell\neq p$ and index $\ell^2$.
Then $\Delta$ is decomposable. 
\rm

\begin{proof} 
We may assume $\ell$ is odd, since if $\ell=2$ the result is a classical theorem of Albert.
Assume the setup \eqref{setup},
let $\alpha=[\Delta]\in\H^2(F,\mu_\ell)$, and let 
$D_\alpha$ be the ramification divisor of $\alpha$ on $X$.
We may assume
$D_\alpha\cup C$ has normal crossings and horizontal components contained in $\D_\S$,
and that we have an element
$\pi\in F$ as in \eqref{setup} with $\div(\pi)=C+E$, where $E$ is horizontal
and avoids the nodal points of $D_\alpha\cup C$.
By Grunwald-Wang's theorem
there exist elements $\theta_i\in\H^1(\kappa(C_i),\Z/\ell)$
such that
\begin{enumerate}
\item[a)]
$\partial_z(\theta_i)=0$ when $z\in C_i$ is a singular point of $C\cup D_\alpha\cup E$.
\item[b)]
$\theta_i(z)=\theta_j(z)$ whenever $z\in C_i\cap C_j$.
\item[c)]
$\theta_i(z)=0$ at all $z\in E\cap C_i$.
\item[d)]
If $z\in D_\alpha\cap\S$ then
\newline\indent
(i)
$\br{\partial_{C_i}(\alpha)(z)-\theta_i(z)}=\br{-\theta_i(z)}$ if $z\in C_i$ is
a curve point of $D_\alpha$;
\newline\indent
(ii)
$\theta_i(z)=0$ if $z\in C_i$ is a not-hot nodal point of $D_\alpha$;
\newline\indent
(iii)
$\br{\partial_{C_i}(\alpha)(z)-\theta_i(z)}=\br{\partial_{C_j}(\alpha)(z)-\theta_i(z)}$
if $z\in C_i\cap C_j$ is a hot point of $\alpha$.
\item[e)]
If $z\in D_\alpha\backslash\S$ then
\newline\indent
(i)
$\br{-\theta_i(z)}=\br{\partial_D(\alpha)(z)}$ if $z\in C_i$ is a curve point of $D_\alpha$;
\newline\indent
(ii)
$\theta_i(z)=0$ if $z\in C_i$ is a not-hot nodal point of $D_\alpha$;
\newline\indent
(iii)
$\br{\partial_{C_i}(\alpha)(z)-\theta_i(z)}=\br{\partial_D(\alpha)(z)}$ if $z\in C_i\cap D$ is a hot point of $\alpha$.
\end{enumerate}
Note that we may arrange (a), (b), and (c) since the given point sets are finite;
(d)(i,iii) and (e)(i,iii) makes sense since the given residues of $\alpha$ are unramified at the given $z$;
(d)(i,iii) and (e)(iii) are possible since $\ell$ is odd;
(c) does not conflict with (d)(i,iii) and (e)(i,iii) since $E$ avoids the nodal points of 
$D_\alpha\cup C$; and
(b) does not conflict with (d)(ii,iii) and (e)(ii,iii) by symmetry.

The $\theta_i$ are unramified with equal values at all nodal points $z\in\S$ by (a,b),
so there exists an element $\theta_C\in\H^1(\O_{C,\S},\Z/\ell)$ 
mapping to $(\theta_1,\dots,\theta_m)\in\H^1(\k(C),\Z/\ell)$ by
Theorem~\ref{lambda}(a).
Note that $\theta_C$ is nonzero by (d)(iii) (or (e)(iii)) since $\alpha$ has at least
one hot point by the hot point criterion \cite[Corollary 5.2]{Sa07}, 
and then $\theta_C(z)=\theta_i(z)$ is necessarily nonzero.

Let $\gamma_1=(\pi)\cdot\lambda(\theta_C)$.
Then $\gamma_1$ and $E$ satisfy the hypotheses of Lemma~\ref{lemma} by (c) and the assumptions on $E$,
hence $D_{\gamma_1}\cup C$ has normal crossings and distinguished horizontal components, 
and since $\theta_C$ is ramified at 
all nodal points of $D_{\gamma_1}$, $\gamma_1$ has no hot points, hence it has index $\ell$
by the hot point criterion.
Write
\begin{align*}
D_\alpha&=C'+H\\
D_{\gamma_1}&=C''+H'\\
\end{align*}
where $C',C''\subset C$, and $H,H'\subset\D_\S$ are distinguished horizontal divisors. 
Set 
\[\gamma_2=\alpha-\gamma_1\]
We intend to show that $\gamma_2$ has index $\ell$.
Since $\theta_C$ is unramified at all singular points of $D_\alpha\cup C$ by (a),
$H'$ avoids all of these points by Lemma~\ref{lemma}, hence $H\cap H'=\varnothing$.
Evidently $D_{\gamma_2}\subset C+H+H'$, hence
$D_{\gamma_2}$ has normal crossings on $X$.
Since $D_{\gamma_2}$ has normal crossings, $\gamma_2$ has index $\ell$ if and only if $\gamma_2$
has no hot points on $X$.

For the following analysis, note
the nodal points $\S_{\gamma_2}$ of $D_{\gamma_2}$ are in $H'$, $\S$, and $H$,
and in the latter two cases $\theta_C$ is unramified at $z$, hence has a value
$\theta_C(z)$.

Suppose $z\in\S_{\gamma_2}$ and $z\not\in D_\alpha$.
Then $z$'s status as a point of $D_{\gamma_2}$ (hot, not hot) 
is the same as its status as a point of $D_{\gamma_1}$, which is not hot since $\gamma_1$
has no hot points. 

Suppose $z\in \S_{\gamma_2}\cap D_\alpha\cap H'$.
Then $z\in C_i\cap D$ for some $C_i$ and some prime divisor
$D\subset H'$, and $\theta_C$ is ramified at $z$ by Lemma~\ref{lemma}.
We have $\partial_D(\alpha)=0$ since $D\not\subset D_\alpha$,
and by Lemma~\ref{lemma}
$$
\partial_{z}(\partial_D(\gamma_2))=\partial_{z}(-\partial_D(\gamma_1))
=\partial_{z}(\partial_{z}(\theta_C)\cdot(\pi))
=v_{z}(\bar\pi)\partial_{z}(\theta_C)
$$
where $\bar\pi$ is the image of $\pi$ in $\kappa(D)$.
Since $\div(\pi)$ has normal crossings at $z$, $v_{z}(\bar\pi)=1$,
hence $\partial_{z}(\partial_D(\gamma_2))=\partial_{z}(\theta_C)\neq 0$.
Since $\partial_z(\partial_D(\gamma_2))\neq 0$, $z$ is not a hot point of $\gamma_2$.

Suppose $z\in \S_{\gamma_2}\cap D_\alpha\cap\S$.
If $z$ is a curve point of $D_\alpha$ on $C_i$, i.e., $\partial_{C_j}(\alpha)=0$
where $C_j$ is the other component of $C$ at $z$, then 
$\partial_{C_i}(\gamma_2)(z)=\partial_{C_i}(\alpha)(z)-\theta_C(z)$
and $\partial_{C_j}(\gamma_2)(z)=-\theta_C(z)$ by Lemma~\ref{lemma}, 
and so $z$ is not a hot point of $\gamma_2$ by (d)(i).
If $z$ is a nodal point of $D_\alpha$ on $C_i\cap C_j$,
then $\theta_C(z)=0$ if $z$ is not a hot point of $\alpha$ by (d)(ii),
so that the status of $z$ for $\gamma_2$ is the same as for $\alpha$ (not hot);
otherwise
$\br{\partial_{C_i}(\gamma_2)(z)}=\br{\partial_{C_j}(\gamma_2)(z)}$
by (d)(iii), hence $z$ is not a hot point for $\gamma_2$ in any case.

Suppose $z\in\S_{\gamma_2}\cap D_\alpha\cap H$.
Assume $z\in C_i\cap D$ for a prime divisor $D\subset H$.
Then $\partial_D(\gamma_1)=0$ since $H\cap H'=\varnothing$.
If $z$ is a curve point of $D_\alpha$, i.e., $\partial_{C_i}(\alpha)=0$, 
then $\partial_{C_i}(\gamma_2)(z)=-\theta_C(z)$ and
$\partial_D(\gamma_2)=\partial_D(\alpha)$, 
hence $\br{\partial_{C_i}(\gamma_2)(z)}=\br{\partial_D(\gamma_2)(z)}$ 
by (e)(i), so $z$ is not a hot point for $\gamma_2$.
If $z$ is a not-hot nodal point of $D_\alpha$ then 
$\theta_C(z)=0$ by (e)(ii), so the status of $z$ is unchanged (not hot) for $\gamma_2$.
If $z$ is a hot point of $\alpha$ then
$\br{\partial_{C_i}(\alpha)(z)-\theta_C(z)}=\br{\partial_D(\alpha)(z)}$ by (e)(iii),
hence $z$ is not a hot point for $\gamma_2$.
This completes the analysis.
We conclude $\gamma_2$ has no hot points on $X$,
hence $\gamma_1$ and $\gamma_2$ both have index $\ell$.
 
Let $\Delta_1$ and
$\Delta_2$ be the $F$-division algebras underlying $\gamma_1$ and $\gamma_2$, respectively,
so that $[\Delta]=[\Delta_1\otimes_F\Delta_2]$.
Since $\ind(\Delta_1\otimes_F\Delta_2)=\ind(\Delta)=\ell^2$ and $\ind(\Delta_i)=\ell$,
it follows that $\Delta_1\otimes_F\Delta_2$
is a division algebra, hence $\Delta\isom\Delta_1\otimes_F\Delta_2$.
\end{proof}

As mentioned in \eqref{brauer}, it is known that $\lBrdim(F)=2$ and
that the $\Z/\ell$-length $\lL(F)$ is finite.  
We now have the following.

\Corollary
Let $F$ be a field finitely generated and of transcendence degree one over $\Q_p$,
and let $\ell\neq p$ be a prime.
Then $\lL(F)=2$.
\rm

\begin{proof}
If $\alpha\in\H^2(F,\mu_\ell)$ then the index of $\alpha$ is either $\ell$ or $\ell^2$
by \cite[Theorem 3.4]{Sa97}.  If it is $\ell$, then $\alpha$ is already $\Z/\ell$-cyclic by
\cite[Theorem 5.1]{Sa07}.  If it is $\ell^2$ then $\alpha=\gamma_1+\gamma_2$ 
for classes $\gamma_i$ in $\H^2(F,\mu_\ell)$ of index $\ell$ by Theorem~\ref{theorem}.
These classes are again $\Z/\ell$-cyclic by Saltman's theorem, and the result follows.
\end{proof}

Saltman proved that all $F$-division algebras of prime degree $\ell$ are cyclic crossed products
in \cite{Sa07}, and Suresh proved the prime period case when $F$ contains the $\ell$-th roots of
unity in \cite{Sur10}.
We now have the prime period case in general:

\Corollary
Let $F$ be a field finitely generated and of transcendence degree one over $\Q_p$, and let $\Delta$
be a division algebra of prime period $\ell\neq p$.  Then $\Delta$ is a crossed product.
\rm

\begin{proof}
The index of $\Delta$ is either $\ell$ or $\ell^2$ by \cite[Theorem 3.4]{Sa97}.  If it is $\ell$, then
$\Delta$ is a cyclic crossed product
by \cite[Theorem 5.1]{Sa07}.  If it is $\ell^2$ then $\Delta=\Delta_1\otimes_F\Delta_2$
by Theorem~\ref{theorem}, and each $\Delta_i$ is cyclic by Saltman's theorem.
Let $L_i/F$ be a cyclic Galois maximal subfield of $\Delta_i$.
Then $L=L_1\otimes_F L_2$ is a commutative Galois subalgebra of $\Delta$ of degree $\ell^2$.  
Since $L$ obviously splits $\Delta$,
$\Delta$ is a crossed product by \cite[Theorem 7.2]{Sa99}.
\end{proof}

\bibliographystyle{abbrv} 
\bibliography{hnx.bib}

\begin{thebibliography}{10}

\bibitem{ABGV}
A.~Auel, E.~Brussel, S.~Garibaldi, and U.~Vishne.
\newblock Open problems on central simple algebras.
\newblock {\em Transform. Groups}, 16(1):219--264, March 2011.

\bibitem{BMT}
E.~Brussel, K.~McKinnie, and E.~Tengan.
\newblock Indecomposable and noncrossed product division algebras over function
  fields of smooth $p$-adic curves.
\newblock {\em Adv. in Math.}, 226:4316--4337, 2011.

\bibitem{BT11}
E.~Brussel and E.~Tengan.
\newblock Tame covers and cohomology of relative curves over complete discrete
  valuation rings, with applications to the {B}rauer group.
\newblock http://arxiv.org/abs/1104.0439, 2011.

\bibitem{GMS}
S.~Garibaldi, A.~Merkurjev, and J.-P. Serre.
\newblock {\em Cohomological invariants in {G}alois cohomology}, volume~28 of
  {\em University Lecture Series}.
\newblock Amer.\ Math.\ Soc., 2003.

\bibitem{Liu}
Q.~Liu.
\newblock {\em Algebraic Geometry and Arithmetic Curves}, volume~6 of {\em
  Oxford Graduate Texts in Mathematics}.
\newblock Oxford University Press, Oxford, 2002.
\newblock Translated from the French by Reinie Ern{\'e}, Oxford Science
  Publications.

\bibitem{Mer83}
A.~Merkurjev.
\newblock Brauer groups of fields.
\newblock {\em Comm. Alg.}, 11:2611--2624, 1983.

\bibitem{MS83}
A.~Merkurjev and A.~Suslin.
\newblock ${K}$-cohomology of {S}everi-{B}rauer varieties and the norm residue
  homomorphism.
\newblock {\em Math. USSR Izv.}, 21(2):307--340, 1983.

\bibitem{Sa97}
D.~Saltman.
\newblock Division algebras over $p$-adic curves.
\newblock {\em J. Ramanujan Math. Soc.}, 12:25--47, 1997.
\newblock see also the erratum \cite{Sa98} and survey \cite{Br10}.

\bibitem{Sa98}
D.~Saltman.
\newblock Correction to division algebras over $p$-adic curves.
\newblock {\em J. Ramanujan Math. Soc.}, 13:125--129, 1998.

\bibitem{Sa99}
D.~Saltman.
\newblock {\em Lectures on {D}ivision {A}lgebras}, volume~94 of {\em CBMS
  Regional Conference Series in Mathematics}.
\newblock American Mathematical Society, Providence, RI, 1999.

\bibitem{Sa07}
D.~Saltman.
\newblock Cyclic algebras over $p$-adic curves.
\newblock {\em J. Algebra}, 314:817--843, 2007.

\bibitem{Sur10}
V.~Suresh.
\newblock Bounding the symbol length in the {G}alois cohomology of function
  fields of {$p$}-adic curves.
\newblock {\em Comment. Math. Helv.}, 85(2):337--346, 2010.

\end{thebibliography}

\end{document}